\documentclass[11pt,leqno]{article}

\usepackage[english]{babel}

\usepackage{amsmath}

\usepackage{amssymb}

\usepackage{xypic}

\usepackage[all]{xy}

\usepackage{latexsym}

\usepackage{amsfonts}

\usepackage{epsf}

\usepackage[dvips]{graphicx}

\usepackage[latin1]{inputenc}

\usepackage[T1]{fontenc}

\numberwithin{equation}{section}

\newcommand{\C}{\mathcal{C}}

\newcommand{\D}{\mathcal{D}}

\newcommand{\M}{\mathcal{M}}

\renewcommand{\mod}{\mathrm{Mod}}

\newcommand{\op}{\mathrm{Op}}

\newcommand{\rc}{\mathbb{R}\textrm{-}\mathrm{c}}

\newcommand{\CC}{\mathbb{C}}

\newcommand{\R}{\mathbb{R}}

\newcommand{\Z}{\mathbb{Z}}

\newcommand{\OO}{\mathcal{O}}

\newcommand{\E}{\mathcal{E}}

\newcommand{\ltens}{\overset{L}{\otimes}}

\newcommand{\iso}{\stackrel{\sim}{\to}}

\newcommand{\dbtxr}{\mathcal{D}\mathit{b}^t_{X_\mathbb{R}}}

\newcommand{\dbt}{\mathcal{D}\mathit{b}^t}

\newcommand{\db}{\mathcal{D}\mathit{b}}

\newcommand{\ot}{\mathcal{O}^t}

\newcommand{\ddxy}{\mathcal{D}_{\xmenoy}}

\newcommand{\ddyx}{\mathcal{D}_{\ymenox}}

\newcommand{\xmenoy}{{X}\rightarrow{Y}}

\newcommand{\ymenox}{{Y}\leftarrow{X}}

\newcommand{\muh}{\mu\mathit{hom}}

\newcommand{\OWX}{\OO^\mathrm{w}_X}
\newcommand{\OWY}{\OO^\mathrm{w}_Y}
\newcommand{\CWXR}{\C^{{\infty ,\mathrm{w}}}_{X_\R}}

\newcommand{\CWM}{\C^{{\infty ,\mathrm{w}}}_M}

\newcommand{\wtens}{\overset{\mathrm{w}}{\otimes}}

\newcommand{\ol}{\OO^\lambda}

\newcommand{\vchar}{\mathrm{Char}}

\newcommand{\sol}{\mathcal{S}ol}

\newcommand{\rh}{\mathit{R}\mathcal{H}\mathit{om}}

\newcommand{\ho}{\mathcal{H}\mathit{om}}

\newcommand{\Ho}{\mathrm{Hom}}

\renewcommand{\dim}{\textbf{Proof.}}

\newcommand{\qed}{\nopagebreak \phantom{} \hfill $\Box$ \\}

\newcommand{\ppi}{\dot{\pi}}

\newcommand{\supp}{\mathrm{supp}}

\newcommand{\RP}{\mathbb{R}^{{\scriptscriptstyle{+}}}}

\newcommand{\imin}[1]{#1^{-1}}

\newcommand{\lind}[1]{\underset{#1}{\underrightarrow{\lim}}}

\newcommand{\dt}[3]{{#1} \to {#2} \to {#3} \stackrel{+}{\to}}

\newtheorem{teo}{Theorem}[subsection]

\newtheorem{df}[teo]{Definition}

\newtheorem{cor}[teo]{Corollary}

\newtheorem{oss}[teo]{Remark}

\newtheorem{prop}[teo]{Proposition}

\newtheorem{lem}[teo]{Lemma}

\author{Luca Prelli}
\title{\bf{CAUCHY-KOWALESKAYA-KASHIWARA THEOREM WITH GROWTH CONDITIONS}}
\date{}

\usepackage{latexsym}
\begin{document}

\maketitle

\begin{abstract}
Here we prove the Cauchy-Kowaleskaya-Kashiwara theorem for holomorphic functions with growth conditions.
\end{abstract}

\tableofcontents

\addcontentsline{toc}{section}{\textbf{Introduction}}

\section*{Introduction}

Let $X$ be a complex manifold. In \cite{KS01} the authors proved that some functional spaces which are not defined by local properties, as tempered and Whitney holomorphic functions are objects of the derived category of sheaves on the subanalytic site associated to $X$, i.e. the site whose objects are subanalytic open subsets of $X$ and the coverings are locally finite. Moreover, they have a $\rho_!\D_X$-module structure (locally, a section of $\Gamma(U;\rho_!\D_X)$ may be written as $P=\sum_{|\alpha| \leq m}a_\alpha(z)\partial^\alpha_z$ with $a_\alpha(z)$ holomorphic on $\overline{U}$).

It is natural to ask if classical results for $\D$-modules remain true in this framework. Here we prove the Cauchy-Kowaleskaya-Kashiwara theorem for holomorphic functions with growth conditions, i.e. we prove that given a coherent $\D_Y$-module $\M$, and a morphism of complex manifolds $f:X \to Y$  which is non characteristic for $\M$, then
$$
\imin f \rh_{\rho_!\D_Y}(\rho_!\M,\OO^\lambda_Y) \simeq \rh_{\rho_!\D_X}(\rho_!\imin {\underline{f}} \M,\ol_X),
$$
where $\lambda$ denote growth conditions on the sheaf of holomorphic functions. The idea of the proof is the following: we divide the proof in two parts. In the first one we prove that the characteristic variety of a coherent $\D_Y$-module $\M$ coincide with the microsupport of $\rh_{\rho_!\D_Y}(\rho_!\M,\OO^\lambda_Y)$. Then, since $f$ is non characteristic for $\M$, we have the isomorphism $\imin f \simeq f^![d]$, where $d$ denotes the difference of the complex dimensions of $X$ and $Y$.
In the second part we use inverse image formulas for $\ol_Y$ to finish the proof of the theorem.\\

In more details, the contents of this paper are as follows.

In \textbf{Section \ref{1}} we recall the definition of the subanalytic site and the construction of microlocalization of subanalytic sheaves.

In \textbf{Section \ref{2}} we recall some results on $\D$-modules and some example of subanalytic sheaves of differential operators and their microlocalization.

The first part of the proof of the theorem is the aim of \textbf{Section \ref{3}}. First we prove that the support of the microlocalization of a subanalytic sheaf $F$ is contained in its microsupport and we use this fact to prove that $\imin f F \simeq f^!F \otimes \omega_{X|Y}$ if $f$ is non characteristic for $SS(F)$. Then we prove that
the characteristic variety of a coherent $\D_Y$-module $\M$ coincide with the microsupport of the sheaf of solutions of $\M$ in $\ol_Y$.

Using these results we are ready to prove the theorem in \textbf{Section \ref{4}}, where we show the second step of the proof using inverse image formulas for $\ol_Y$.\\

\noindent \textbf{Acknowledgments.} We thank Pierre Schapira for his valuable suggestions on many aspects of this subject. We thank Andrea D'Agnolo for his remarks and comments during the preparation of this work. We thank Teresa Monteiro Fernandes for the useful discussions we had at CAUL.

\section{Subanalytic sheaves}\label{1}

In this section we recall the definition of the subanalytic site and the construction of microlocalization of subanalytic sheaves. References are made to \cite{KS90} and \cite{Pr1} for the theory of sheaves on subanalytic sites, and to \cite{Pr07} for the microlocalization of subanalytic sheaves.

\subsection{Sheaves on subanalytic sites}\label{1.1}

Let $X$ be a real analytic manifold and let $k$ be a field. Denote
by $\op(X_{sa})$ the category of subanalytic subsets of $X$. One
endows $\op(X_{sa})$ with the following topology: $S \subset
\op(X_{sa})$ is a covering of $U \in \op(X_{sa})$ if for any
compact $K$ of $X$ there exists a finite subset $S_0\subset S$
such that $K \cap \bigcup_{V \in S_0}V=K \cap U$. We will call
$X_{sa}$ the subanalytic site.\\

Let $\mod(k_{X_{sa}})$ denote the category of sheaves on $X_{sa}$
and let $\mod_{\rc}(k_X)$  be the abelian category of
$\R$-constructible sheaves on $X$. We denote by $\rho: X \to X_{sa}$ the natural morphism of sites.
We have functors
\begin{equation*}
\xymatrix{\mod_{\rc}(k_X) \subset \mod(k_X)   \ar@ <2pt>
[r]^{\hspace{1cm}\rho_*} &
  \mod(k_{X_{sa}}) \ar@ <2pt> [l]^{\hspace{1cm}\imin \rho}. }
\end{equation*}
The functor $\imin \rho$ admits a left adjoint, denoted by $\rho_!$. The sheaf
$\rho_!F$ is the sheaf associated to the presheaf $\op(X_{sa}) \ni
U \mapsto F(\overline{U})$.

The functor $\rho_*$ is fully faithful and exact on
$\mod_{\rc}(k_X)$ and we identify $\mod_{\rc}(k_X)$ with its
image in $\mod(k_{X_{sa}})$ by $\rho_*$.\\

\begin{teo} Let $F \in \mod(k_{X_{sa}})$. Then there exists a filtrant inductive system $\{F_i\}$ in $\mod_{\rc}(k_X)$ such that $F \simeq \lind i \rho_* F_i$.
\end{teo}

Let $X,Y$ be two real analytic manifolds, and let $f:X \to Y$ be a
real analytic map. The functors $\ho$, $\otimes$, $\imin f$ and $f_*$ are always
defined for sheaves on Grothendieck topologies. For subanalytic
sheaves we can also define the functor of proper direct image
$f_{!!}$. The functor $Rf_{!!}$ admits a right adjoint, denoted by
$f^!$, and we get the usual isomorphisms like projection formula and
base change formula.\\

\subsection{Microlocalization of subanalytic sheaves}\label{1.2}

Let $E$ be a vector bundle over a real analytic manifold $Z$
endowed with the natural action $\mu$ of $\RP$, the multiplication
on the fibers. Let $U$ be an open subset of $E$. We say that $U$
is $\RP$-connected if its intersections with the orbits of $\mu$
are connected. We denote $\RP U$ the conic open set associated to
$U$ (i.e. $\RP U=\mu(U,\RP)).$\\

\begin{df} A sheaf $F$ on $E_{sa}$ is
said conic if $\Gamma(\RP U;F) \iso \Gamma(U;F)$ for each
$\RP$-connected relatively compact open subanalytic subset $U$ of
$E$. We call $D^b_{\RP}(k_{E_{sa}})$ the
subcategory of $D^b(k_{E_{sa}})$ consisting of objects with conic
cohomology.\\
\end{df}

Let $E^*$ be the dual vector bundle and  consider the projections
$p_1,p_2$ from $E \times_Z E^*$ to $E$ and $E^*$ respectively. Let
$P':=\{(x,y) \in E \times_Z E^*;\; \langle x,y\rangle \leq 0\}$.
As in classical sheaf theory, one can define the Fourier-Sato transform
and the inverse Fourier-Sato transform
\begin{eqnarray*}
(\cdot)^{\land}:D^b_{\RP}(k_{E_{sa}}) \to D^b_{\RP}(k_{E_{sa}^*}),
&&
F^{\land}=Rp_{2!!}(p_1^{-1}F)_{P'},\\
(\cdot)^{\vee}:D^b_{\RP}(k_{E_{sa}^*})\to D^b_{\RP}(k_{E_{sa}}),&&
F^{\vee}=Rp_{1*}\mathrm{R}\Gamma_{P'} p_2^!F.
\end{eqnarray*}
The functors ${}^{\wedge}$ and
${}^{\vee}$ are equivalence of categories, inverse to each others.\\

Let $X$ be a real analytic manifold and let $M$ be a closed
submanifold of $X$. We denote by $T_MX
\overset{\tau}{\to} M$ the normal bundle and by $T^*_MX
\overset{\pi}{\to} M$ the conormal bundle.\\

We consider the normal deformation of $X$, i.e. an analytic
manifold $\widetilde{X}_M$, an application $(p,t):\widetilde{X}_M
\to X \times \R$, and an action of $\R \setminus \{0\}$ on
$\widetilde{X}_M$ $(\widetilde{x},r) \mapsto \widetilde{x} \cdot
r$ such that $\imin p (X \setminus M)\simeq(X \setminus M)
     \times (\R \setminus \{0\})$,
  $\imin t (c)\simeq X$ for each $c \neq 0$ and $\imin t (0)\simeq
  T_MX$.
Let $s:T_MX \hookrightarrow \widetilde{X}_M$ be the inclusion,
$\Omega$ the open subset of $\widetilde{X}_M$ defined by
$\{t>0\}$, $i_\Omega: \Omega \hookrightarrow \widetilde{X}_M$ and
$\widetilde{p}=p \circ i_\Omega$. We get a commutative diagram
\begin{equation*}
\xymatrix{T_MX \ar[r]^s \ar[d]^\tau & \widetilde{X}_M \ar[d]^p &
\Omega \ar[l]_{\ \ i_\Omega} \ar[dl]^{\widetilde{p}} \\
M \ar[r]^{i_M} & X. & }
\end{equation*}

\begin{df}  The specialization along $M$ is the functor
\begin{eqnarray*}
\nu^{sa}_M:D^b(k_{X_{sa}}) & \to & D^b_{\RP}(k_{T_MX_{sa}}) \\
\nu^{sa}_M(F) & = & \imin sR\Gamma_\Omega\imin p F.
\end{eqnarray*}
 The microlocalization along $M$ is the Fourier-Sato
transform of the specialization, i.e.
\begin{eqnarray*}
\mu^{sa}_M:D^b(k_{X_{sa}}) & \to & D^b_{\RP}(k_{T^*_MX_{sa}}) \\
\mu^{sa}_M F & = & (\nu^{sa}_M F)^\wedge.
\end{eqnarray*}
\end{df}
These definitions are compatible with
the classical definitions of specialization and microlocalization of \cite{KS90}: we have the
isomorphism $\imin \rho \circ \nu^{sa}_M \circ R\rho_* \simeq \nu_M$ and $\imin \rho \circ \mu^{sa}_M \circ R\rho_* \simeq \mu_M$.

We have the Sato's triangle for subanalytic sheaves:
$$
F|_M \otimes \omega_{M/X} \to {\rm R}\Gamma_MF|_M \to
R\dot{\pi}_*\mu^{sa}_MF \stackrel{+}{\to}
$$
where $\dot{\pi}$ is the restriction of $\pi$ to $T^*_MX \setminus
M$. \\

Let $\Delta$ be the diagonal of $X \times X$, and denote by
$\delta$ the diagonal embedding. The normal deformation of the
diagonal in $X \times X$ can be visualized by the following
diagram
\begin{equation}\label{normdef}
\xymatrix{TX \ar[r]^{\hspace{-0.8cm}\sim} & T_{\Delta}(X \times X)
\ar[r]^{\hspace{0.3cm}s} \ar[d]^{\tau_X} & \widetilde{X \times X}
\ar[d]^p & \Omega
\ar[l]_{\hspace{0.5cm}i_\Omega} \ar[dl]^{\widetilde{p}} \\
& \Delta \ar[r]^{\delta} \ar[dr]^{\sim} & X \times X \ar@ <2pt>
[d]^{q_2} \ar@ <-2pt> [d]_{q_1} & \\
& & X. & }
\end{equation}

\begin{df} Let $F,G \in D^b(k_{X_{sa}})$. We set
$$\muh^{sa}(F,G):=\mu^{sa}_\Delta \rh(\imin q_2 F, q_1^!G).$$
\end{df}

\section{$\D$-modules and $\E$-modules}\label{2}

In this section we recall some example of subanalytic sheaves of differential operators and their microlocalization. Reference are made to \cite{Ka03} and \cite{Sc85} for an introduction to $\D$-modules and $\E$-modules respectively. Tempered microlocalization has been studied in detail in \cite{An94} and we refer to \cite{Co98} for the definition of formal microlocalization. The link between tempered and formal microlocalization and subanalytic sheaves can be found in \cite{Pr07}.

\subsection{Notations and review}\label{2.1}

Let $k=\CC$ and let $X$ be a complex analytic manifold. We denote by $\D_X$ the sheaf of rings of differential operators. Locally, a section of $\Gamma(U;\D_X)$ may be written as $P=\sum_{|\alpha| \leq m}a_\alpha(z)\partial^\alpha_z$ with $a_\alpha(z)$ holomorphic on $U$. We denote by $\mod(\D_X)$ the category of sheaves of $\D_X$-modules.\\


Let $f:X \to Y$ be a morphism of complex analytic manifolds and let $\ddxy=\OO_X \otimes_{\imin f\D_Y}\D_Y$ be the transfer bimodule of $f$. The inverse image of a $\D_Y$-module $\M$ is defined by
$$
\imin{\underline{f}}\M=\ddxy \otimes_{\imin f \D_Y} \imin f \M.
$$

Let $T^*X \stackrel{\pi}{\to} X$ be the cotangent bundle. We denote by $\E_X$ the sheaf of rings of microdifferential operators.

\begin{df} The characteristic variety $\vchar(\M)$ of a $\D_X$-module $\M$ is the support of $\E_X \otimes_{\imin \pi \D_X} \imin \pi \M$.
\end{df}


Let $f:X \to Y$ be a morphism of complex analytic manifolds and let $f_\pi:X \times_Y T^*Y \to T^*Y$ be the base change map.

\begin{df} Let $f:X \to Y$ be a morphism of complex manifold. Then $f$ is non characteristic for $\M$ if
$$
\imin f_\pi(\vchar(\M)) \cap T^*_XY \subseteq X \times_Y T^*_YY.
$$
\end{df}


We recall the following result

\begin{teo} Let $f:X \to Y$ be a morphism of complex analytic manifolds and let $\M$ be a coherent $\D_Y$-module. Assume that $f$ is non characteristic for $\M$. Then $\imin {\underline{f}} \M$ is a coherent $\D_X$-module and there is an isomorphism
$$
\imin f \rh_{\D_Y}(\M,\OO_Y) \simeq \rh_{\D_X}(\imin {\underline{f}}\M,\OO_X).
$$
\end{teo}

\subsection{Microlocalization with growth conditions}\label{2.2}

 Let $M$ be a real analytic
manifold. One denotes by $\db_M$ and $\C^\infty_M$ the sheaves of
Schwartz's distributions and $\C^\infty$ functions respectively.
We recall the definitions of the sheaves of tempered distributions
$\dbt_M$ and Whitney $\C^\infty$ functions $\CWM$ on $M_{sa}$. We have:
\begin{eqnarray*}
\Gamma(U;\dbt_M) & = & \Gamma(M;\db_M)/\Gamma_{M\setminus
U}(M;\db_M), \\
\Gamma(U;\CWM) & = &
\Gamma(M;\C^\infty_M)/\Gamma(M;\mathcal{I}^\infty_{M\setminus
  U}),
\end{eqnarray*}
where $U$ is a locally cohomologically trivial subanalytic subset
and $\Gamma(M;\mathcal{I}^\infty_{M\setminus
  U})$ denotes the space of $\C^\infty$ functions vanishing on
$M\setminus U$ with infinite order.\\

Now let $X$ be a complex manifold, $X_\R$ the underlying real
analytic manifold and $\overline{X}$ the complex conjugate
manifold. One denotes by $\OO_X$ the sheaf of holomorphic functions on $X$ and by $\D_X$ the sheaf of finite order
differential operators with holomorphic coefficients. We denote by $\OO^\lambda_X$, $\lambda=t,\mathrm{w},\omega$ the objects of $D^b(\rho_!\D_X)$ defined by:
\begin{eqnarray*}
\OO^t_X & = & \rh_{\rho_!\D_X}(\rho_!\OO_{\overline{X}},\dbtxr), \\
\OWX
& = & \rh_{\rho_!\D_X}(\rho_!\OO_{\overline{X}},\CWXR), \\
\OO^\omega_X & = & \rho_!\OO_X.
\end{eqnarray*}
We shall need the following isomorphism of \cite{KS01}. Let $F \in D^b_{\rc}(\CC_X)$, then
\begin{equation}\label{homega}
\imin \rho \rh(F,\OO^\omega_X) \simeq D'F \otimes \OO_X,
\end{equation}
where $D'(F)=\rh(F,\CC_X)$.\\

Let us consider the normal deformation of the diagonal in $X
\times X$ as in diagram \eqref{normdef}. Microlocalization of tempered and Whitney holomorphic functions correspond to the functors of tempered and formal microlocalization. In fact, let $F \in
D^b_{\rc}(\CC_X)$, we have the isomorphisms
\begin{eqnarray*}
\imin \rho \muh^{sa}(F,\ot_X) & \simeq & t\muh(F,\OO_X), \\
\imin \rho \muh^{sa}(F,\OWX) & \simeq & (D'F
\underset{\mu}{\wtens} \OO_X)^a,
\end{eqnarray*}
where $t\muh$ and $\underset{\mu}{\wtens}$ are the functors of tempered and formal microlocalization (see \cite{An94} and \cite{Co98} for details) and $(\cdot)^a$ denotes
the direct image for the antipodal map. We shall need the following result

\begin{teo}\label{muhEmod} Let $F \in D^b_{\rc}(\CC_X)$. Then $H^k\muh(F,\OO^\lambda_X)$ has a structure of $\E_X$-module, for $\lambda=\varnothing,t,\mathrm{w}$ and for any $k \in \Z$.
\end{teo}

\section{Microsupport and characteristic variety}\label{3}

In this section we expose some results on microsupport of subanalytic sheaves and its relation with the functor of microlocalization and the characteristic variety of a $\D$-module. References are made to \cite{KS03} for the construction of the microsupport (for ind-sheaves), and to \cite{Ma08} for the functorial properties of the microsupport.

\subsection{Microsupport of subanalytic sheaves}\label{3.1}

Let $X$ be a real analytic manifold and let $T^*X \stackrel{\pi}{\to} X$ be the cotangent bundle. We recall the following two equivalent definitions of microsupport of a subanalytic sheaves of \cite{KS03}. For the notion of microsupport for classical sheaves we refer to \cite{KS90}.

\begin{df} The microsupport of $F \in D^b(k_{X_{sa}})$, denoted by $SS(F)$ is the subset of $T^*X$ defined as follows. Let $p \in T^*X$, then $p \notin SS(F)$ if one of the following equivalent conditions is satisfied.
\begin{itemize}
\item[(i)] There exists a conic neighborhood $U$ of $p$ and a small filtrant system $\{F_i\}$ in $C^{[a,b]}(\mod_{\rc}(k_X))$ with  $SS(F_i) \cap U=\varnothing$ such that $F$ is quasi-isomorphic to $\lind i \rho_* F_i$ in a neighborhood of $\pi(p)$.
\item[(ii)] There exists a conic neighborhood $U$ of $p$ such that for any $G \in D_{\rc}^b(k_X)$ with $\supp(G)\subset\subset \pi(U)$ and such that $SS(G) \subset U \cup T^*_XX$, one has $\Ho_{D^b(k_{X_{sa}})}(G,F)=0$.
\end{itemize}
\end{df}

\begin{oss} In \cite{KS03} microsupport was defined for ind-sheaves. The above definition follows from the equivalence between subanalytic sheaves and ind-$\R$-constructible sheaves (see \cite{KS01} for details).
\end{oss}


Let $M$ be a real closed submanifold of $X$.

\begin{prop}\label{muSS(F)} Let $F \in D^b(k_{X_{sa}})$. Then $\supp (\mu_M^{sa} F) \subseteq SS(F) \cap T^*_MX$.
\end{prop}
\dim\ \ Let $F \in D^b(k_{X_{sa}})$ and let $p \notin SS(F)$. There exists conic neighborhood $U$ of $p$ and a small filtrant system $\{F_i\}$ in $C^{[a,b]}(\mod_{\rc}(k_X))$ with  $SS(F_i) \cap \overline{U}=\varnothing$ such that there exists $W \in \op(X_{sa})$ with $U \subseteq \imin \pi (W)$ and $F_W \simeq \lind i \rho_* F_i$.
We have $H^k\mu_M^{sa} F_W \simeq \lind i \rho_* H^k \mu_M F_{iW}$, hence $(\mu_M^{sa} F)|_U=0$ since $\supp (\mu_M F_i) \subseteq SS(F_i)$. \\ \qed

\begin{cor}\label{muhSS(F)G} Let $G \in D^b_{\rc}(k_X)$, $F \in D^b(k_{X_{sa}})$. Then $\supp (\muh^{sa}(F,G)) \subseteq SS(F) \cap SS(G)$.
\end{cor}
The result follows from Proposition \ref{muSS(F)} and the following result of \cite{Ma08}:
$$
SS(\rh(\imin q_1 G,q_2^!F)) \subseteq SS(G)^a \times SS(F).
$$
\qed

Let $f:X \to Y$ be a morphism of real analytic manifolds and denote by $f_\pi: X \times_Y T^*Y \to T*Y$ the base change map.

\begin{df} Let $f:X \to Y$ be a morphism of real analytic manifolds and let $F \in D^b(k_{Y_{sa}})$. One says that $f$ is non characteristic for $SS(F)$ if
$$
\imin f_\pi(SS(F)) \cap T^*_XY \subseteq X \times_Y T^*_YY.
$$
If $f$ is a closed embedding $X$ is said to be non characteristic.
\end{df}

\begin{prop}\label{noncharSS(F)} Let $f:X \to Y$ be a morphism of real analytic manifolds and let $F \in D^b(k_{Y_{sa}})$. Assume that $f$ is non characteristic for $SS(F)$. Then the natural morphism
$$
\imin f F \otimes \omega_{X|Y} \to f^!F
$$
is an isomorphism.
\end{prop}
\dim\ \ We may reduce to the case $f$ closed embedding, hence we have to prove the isomorphism $F|_X \otimes \omega_{X|Y} \simeq R\Gamma_XF|_X$ when $SS(F) \cap T^*_XY \subseteq T^*_YY$. Consider the Sato's triangle
$$
\dt{F|_X \otimes \omega_{X|Y}}{R\Gamma_XF|_X}{R\dot{\pi}_*\mu^{sa}_XF}.
$$
Since $SS(F) \cap T^*_XY \subseteq T^*_YY$ we have $R\dot{\pi}\mu^{sa}_XF=0$ by Proposition \ref{muSS(F)} and the result follows. \\ \qed

\subsection{Characteristic variety}\label{3.2}

Now let us study some applications of the preceding results to
$\D$-modules. We first need the following lemma.

\begin{lem} Let $F \in D^b_{\rc}(\CC_X)$ and let $G \in
D^b(k_{X_{sa}})$ then
\begin{equation}\label{D'Rc}
\imin\rho R\pi_{!!}\muh^{sa}(F,G) \simeq D'F \otimes \imin \rho G.
\end{equation}
\end{lem}
\dim\ \ (i) Let us prove first the isomorphism
\begin{equation}\label{rhD'}
D'F \boxtimes \imin \rho G \iso \imin \rho \rh(\imin {q_1}F,\imin {q_2}G).
\end{equation}
We may reduce to the case $F=\CC_U$ with $D'\CC_U \simeq \CC_{\overline{U}}$. Hence, given $V,W \in \op^c(X_{sa})$, it is enough to prove the isomorphism $H^k(V \times W;(\imin {q_2}G)_{\imin q_1(\overline{U})}) \simeq H^k(V \times W; \mathrm{R}\Gamma_{\imin {q_1}(U)}\imin {q_2}G)$ for each $k \in \Z$. We have
\begin{eqnarray*}
H^k(V \times W;(\imin {q_2}G)_{\imin q_1(\overline{U})}) & \simeq & \lind {U' \supset \overline{U}}H^k(V \cap U' \times W;\imin {q_2}G) \simeq H^k(W;G) \\
H^k(V \times W; \mathrm{R}\Gamma_{\imin {q_1}(U)}\imin {q_2}G) & \simeq & H^k(V \cap U \times W;\imin {q_2}G) \simeq H^k(W;G),
\end{eqnarray*}
where $U' \in \op(X_{sa})$.

(ii) We have
\begin{eqnarray*}
R\pi_{!!}\muh^{sa}(F,G)
& \simeq & \rh(\imin {q_1}F,q_2^!G)|_\Delta \otimes \omega_{\Delta|X
\times X} \\
& \simeq & \rh(\imin {q_1}F,\imin {q_2}G)|_\Delta  \otimes \omega_{X \times X|X}|_\Delta \otimes
\omega_{\Delta|X \times X} \\
& \simeq & \rh(\imin q_1F,\imin {q_2}G)|_\Delta.
\end{eqnarray*}

Applying $\imin \rho$ and \eqref{rhD'} we obtain
\begin{eqnarray*}
\imin \rho R\pi_{!!} \muh^{sa}(F,G) & \simeq & \imin \rho \rh(\imin q_1F,\imin {q_2}G)|_\Delta \\
& \simeq & (D'F \boxtimes  \imin {q_2}\imin \rho G)|_\Delta \\
& \simeq & D'F \otimes \imin \rho G.
\end{eqnarray*}
In the second isomorphism we used the fact that $\imin \rho$
commutes with the functor of inverse image.\\
\qed

 We recall
the notion of elliptic pair of \cite{SS94}. Let $\M$ be a coherent
$\D$-module and $F \in D^b_{\rc}(\CC_X)$, then $(F,\M)$ is an
elliptic pair if
$$ SS (F) \cap \vchar(\M) \subseteq T^*_XX.$$
We consider the sheaf $\OO^\lambda_X$, for
$\lambda=\varnothing,t,\mathrm{w},\omega$.

\begin{prop}\label{ellpair} Let $(F,\M)$ be an elliptic pair. Then we have the
isomorphism
$$
 \rh_{\D_X}(\M,D'F \otimes \OO_X) \iso
 \rh_{\D_X}(\M,\imin \rho \rh(F,\OO^\lambda_X)).
$$
\end{prop}
\dim\ \ If $\lambda=\omega$ the result follows from \eqref{homega}.

Let $\lambda=\varnothing,t,\mathrm{w}$. Let $\delta: \Delta \to X \times X$ be the embedding and
let us consider the Sato's triangle
\begin{eqnarray}\label{Satos}
& \imin \delta\rh(\imin {q_1}F,q_2^!\ol_X) \otimes
\omega_{\Delta|X
\times X}  \to  \delta^!\rh(\imin {q_1}F,q_2^!\ol_X) & \notag \\
 & \to R\dot{\pi}_*\muh^{sa}(F,\ol_X) \stackrel{+}{\to}. &
\end{eqnarray}
We have $\delta^!\rh(\imin {q_1}F,q_2^!\ol_X) \simeq
\rh(F,\ol_X)$. Then applying $\imin \rho$ to the triangle
\eqref{Satos}, using the isomorphism \eqref{D'Rc} and the fact
that $\imin \rho \ol_X \simeq \OO_X$, we obtain
$$
\dt{D'F \otimes \OO_X}{\imin \rho \rh(F,\ol_X)}{\imin \rho
R\dot{\pi}_*\muh^{sa}(F,\ol_X)}.
$$
Applying the functor $\rh_{\D_X}(\M,\cdot)$ we obtain
\begin{eqnarray*}
& \rh_{\D_X}(\M,D'F \otimes \OO_X) \to
 \rh_{\D_X}(\M,\imin \rho \rh(F,\OO^\lambda_X)) & \\
& \to \rh_{\D_X}(\M,\imin \rho R\dot{\pi}_*\muh^{sa}(F,\ol_X))
\stackrel{+}{\to}. &
\end{eqnarray*}
Then it is enough to prove that $\rh_{\D_X}(\M,\imin \rho
R\dot{\pi}_*\muh^{sa}(F,\ol_X))=0$. First remark that since
$R\dot{\pi}_*$ commutes with $\imin \rho$ we have by adjunction
\begin{eqnarray*}
\lefteqn{\rh_{\D_X}(\M,\imin \rho R\ppi_*\muh^{sa}(F,\ol_X))}\\ &
\simeq & R\ppi_*\rh_{\imin \pi \D_X}(\imin \pi \M,\imin \rho
\muh^{sa}(F,\ol_X)).
\end{eqnarray*}
Let $k \in \Z$. By Theorem \ref{muhEmod} we have
\begin{eqnarray*}
\lefteqn{R\ppi_*\rh_{\imin \pi \D_X}(\imin \pi \M,H^k\imin \rho
\muh^{sa}(F,\ol_X))} \\
& \simeq & R\ppi_*\rh_{\E_X}(\E_X \otimes_{\imin \pi \D_X} \imin
\pi \M,H^k\imin \rho \muh^{sa}(F,\ol_X)).
\end{eqnarray*}
We have $\supp(\E_X \otimes_{\imin \pi \D_X} \M)=\vchar(\M)$ and
$\supp(H^k\imin \rho \muh^{sa}(F,\ol_X)) \subseteq SS(F)$
for each $k \in \Z$ by Corollary \ref{muhSS(F)G}. Hence
\begin{equation}\label{Hk=0}
R\ppi_*\rh_{\imin \pi \D_X}(\imin \pi \M,H^k\imin \rho
\muh^{sa}(F,\ol_X))=0
\end{equation}
for each $k \in \Z$ since the pair $(F,\M)$ is elliptic. Let us
suppose that the length  of the bounded complex
$\muh^{sa}(F,\ol_X)$ is $n$ and let us argue by induction on the
truncation $\tau^{\leq i}\muh^{sa}(F,\ol_X)$. If $i=0$ the result
follows from \eqref{Hk=0}. Let us consider the distinguish
triangle
$$
\dt{\tau^{\leq n-1}\imin \rho \muh^{sa}(F,\ol)}{\imin \rho
\muh^{sa}(F,\ol)}{H^n\imin \rho \muh^{sa}(F,\ol)}
$$
and apply the functor $R\dot{\pi}_*\rh_{\imin \pi\D_X}(\imin \pi
\M,\cdot)$. The first term becomes zero by the induction hypothesis
and the third one is zero by \eqref{Hk=0}. Hence
$\rh_{\D_X}(\M,\imin \rho R\dot{\pi}_*\muh^{sa}(F,\ol_X))=0$ and
the result follows.\\
\qed

Let $\M$ be a $\D_X$-module and let $\lambda=\varnothing,t,\mathrm{w},\omega$. One sets for short
$$\sol^\lambda(\M):=\rh_{\rho_!\D_X}(\rho_!\M,\ol_X).$$

\begin{cor}\label{ssvchar} Let $\M$ be a coherent $\D_X$-module. Then
$$
SS(\sol^\lambda(\M))=\vchar(\M).$$
\end{cor}
\dim\ \ Recall that $SS(\sol(\M))=\vchar(\M)$.

(i) $\vchar(\M) \subseteq SS(\sol^\lambda(\M))$ follows from the
fact that $\imin \rho \sol^\lambda(\M)=\sol(\M)$ and $SS(\imin
\rho G) \subseteq SS(G)$ for each $G \in D^b(\CC_{X_{sa}})$.

(ii) $\vchar(\M) \supseteq  SS(\sol^\lambda(\M))$. Let $(x,\xi)
\notin \vchar(\M)=SS(\sol(\M))$ and let $U$ be an open
neighborhood of $x$ such that $U \cap \pi(\vchar(\M))=\varnothing$
and such that for each $F \in D^b_{\rc}(\CC_X)$ with $\supp(F)
\subset \subset U$ one has $\Ho_{D^b(\CC_X)}(F,\sol(\M))=0$. By
Proposition \ref{ellpair} the complexes
\begin{eqnarray*}
\rh(F,\rh_{\D_X}(\M,\OO_X)) & \simeq & \imin \rho \rh(F,
\rh_{\rho_!\D_X}(\rho_!\M,R\rho_*\OO_X))
\\
& \simeq & \imin \rho \rh(F,\rh_{\rho_!\D_X}(\rho_!\M,\ol_X))
\end{eqnarray*}
are all quasi-isomorphic for
$\lambda=\varnothing,t,\mathrm{w},\omega$.
Hence $\Ho_{D^b(\CC_{X_{sa}})}(F,\sol^\lambda(\M))=0$ and $(x,\xi) \notin SS(\sol^\lambda(\M))$.\\
\qed

\section{Cauchy-Kowaleskaya-Kashiwara theorem}\label{4}

In this section we apply the results of the previous one to prove the Cauchy-Kowaleskaya-Kashiwara theorem for holomorphic functions with growth conditions $\lambda=\varnothing,t,\mathrm{w},\omega$. We refer to \cite{Ka03} for the statement and proof of the Cauchy-Kowaleskaya-Kashiwara theorem for holomorphic functions.

\subsection{Statement of the theorem}\label{4.1}

Let $f:X \to Y$ be a morphism of complex manifolds. Set $d=\mathrm{dim}X-\mathrm{dim}Y$. We recall the following isomorphisms of \cite{KS01} and \cite{Pr07}
\begin{eqnarray}
f^!\ot_Y & \simeq & \rho_!\ddyx \underset{\rho_!\D_X}{\ltens}\ot_X[d], \label{iminvot} \\
f^!\OWY & \simeq & \rh_{\rho_!\D_X}(\rho_!\ddxy,\OWX)[2d]. \label{iminvow}
\end{eqnarray}

Our main theorem is the following

\begin{teo}\label{ckk} Let $\M$ be a coherent $\D_Y$-module, and suppose that $f$ is non characteristic for $\M$. Then we have the following isomorphism for $\lambda=\varnothing,t,\mathrm{w},\omega$
$$
\imin f \rh_{\rho_!\D_Y}(\rho_!\M,\OO^\lambda_Y) \simeq \rh_{\rho_!\D_X}(\rho_!\imin {\underline{f}} \M,\ol_X).
$$
\end{teo}

\subsection{Proof of the theorem}\label{4.2}

\begin{prop}\label{division} Let $\M$ be a coherent $\D_Y$-module, and suppose that $f$ is non characteristic for $\M$.  Then we have the following isomorphism for $\lambda=\varnothing,t,\mathrm{w},\omega$
$$
f^! \rh_{\rho_!\D_Y}(\rho_!\M,\ol_Y) \simeq \rh_{\rho_!\D_X}(\rho_!\imin {\underline{f}} \M,\ol_X)[2d].
$$
\end{prop}
\dim\ \ (i) Let $\lambda=t$. Recall that if $\M$ is a coherent $\D_Y$-module and $f$ is non characteristic, then $\imin {\underline{f}} \M$ is a coherent $\D_X$-module and
$$
\imin {\underline{f}} \rh_{\D_Y}(\M,\D_Y) \simeq \rh_{\D_X}(\imin {\underline{f}} \M,\D_X)[d].
$$
We have the chain of isomorphisms
\begin{eqnarray*}
\rh_{\D_X}(\rho_!\imin {\underline{f}}\M,\ot_X)[2d] & \simeq &  \rho_!\rh_{\D_X}(\imin {\underline{f}}\M,\D_X) \otimes_{\rho_!\D_X} \ot_X[2d] \\
& \simeq & \rho_!\imin {\underline{f}} \rh_{\D_Y}(\M,\D_Y) \otimes_{\rho_!\D_X} \ot_X[d] \\
& \simeq & \rho_!\imin f \rh_{\D_Y}(\M,\D_Y) \otimes_{\rho_!\imin f\D_Y} f^!\ot_Y \\
& \simeq & f^!(\rho_!\rh_{\D_Y}(\M,\D_Y) \otimes_{\rho_!\D_Y} \ot_Y) \\
& \simeq & f^!\rh_{\rho_!\D_Y}(\rho_!\M,\ot_Y),
\end{eqnarray*}
where the first and the last isomorphisms follow from the coherence of $\imin {\underline{f}} \M$ and $\M$, and the third one follows from \eqref{iminvot}.\\

\noindent (ii) Let $\lambda=\mathrm{w}$. We have the chain of isomorphisms
\begin{eqnarray*}
f^!\rh_{\rho_!\D_Y}(\M,\OWY) & \simeq & \rh_{\rho\imin f \D_X}(\rho_!\imin f \M,f^!\ot_Y) \\
& \simeq & \rh_{\rho_!\imin f\D_X}(\rho_!\imin f \M,\rh_{\rho_!\D_X}(\rho_!\ddxy,\OWX))[2d] \\
& \simeq & \rh_{\rho_!\D_X}(\rho_!\imin {\underline{f}} \M,\OWX)[2d],
\end{eqnarray*}
where the second isomorphism follows from \eqref{iminvow}.\\

\noindent (iii) Let $\lambda=\varnothing,\omega$. Since $\M$ is coherent and $f$ is non characteristic the result follows from the isomorphism
$$
f^! \rh_{\D_Y}(\M,\OO_Y) \simeq \rh_{\D_X}(\imin {\underline{f}} \M,\OO_X)[2d].
$$
\qed

\bigskip

\noindent \textbf{Proof of Theorem \ref{ckk}.}\ \ By Corollary \ref{ssvchar} $f$ is non characteristic for $SS(\sol^\lambda(\M))$. Hence by Proposition \ref{noncharSS(F)}
$$
f^! \rh_{\rho_!\D_Y}(\rho_!\M,\ol_Y) \simeq \imin f \rh_{\rho_!\D_Y}(\rho_!\M,\ol_Y) [2d].
$$
Then the result follows from Proposition \ref{division}.\\
\qed


\addcontentsline{toc}{section}{

\end{document}